\documentclass[a4paper]{article}
\usepackage{color}
\usepackage{amsmath, amssymb, theorem,latexsym}
\usepackage{graphicx}
\usepackage{epsf}
\usepackage[english]{babel}
\theoremstyle{plain}
\newtheorem{theorem}{Theorem}[section]

\newtheorem{corollary}[theorem]{Corollary}

\newcommand{\R}{\mathbb{R}}
\newcommand{\D}{\mathbb{D}}

\title{The weak Pleijel theorem with geometric control}
\author{Pierre B\'erard \\ Institut Fourier, Universit\'{e} de Grenoble and
CNRS\\ B.P.74 ~~F 38 402 Saint Martin d'H\`{e}res Cedex, France.\\
\texttt{pierrehberard@gmail.com}\\
and \\
Bernard Helffer \\
Laboratoire de Math\'ematiques Jean Leray\\
CNRS, Universit\'e de Nantes.\\
2 rue de la Houssini\`ere, F44322 Nantes Cedex 3, France\\
Laboratoire de Math\'ematiques d'Orsay, Univ. Paris-Sud, CNRS\\
Universit\'e Paris-Saclay,  France.
\\
\texttt{Bernard.Helffer@univ-nantes.fr}
}

\date{To appear in Journal of Spectral Theory \textbf{6} (2016)}

\begin{document}
\maketitle

\begin{center}
\emph{In memory of Y. Safarov}
\end{center}%
\vspace{1cm}

\begin{abstract}
Let $\Omega\subset \mathbb R^d\,, d\geq 2$, be a bounded open set, and denote by $\lambda_j(\Omega), j\ge 1$, the eigenvalues of the Dirichlet Laplacian arranged in nondecreasing order, with multiplicities. The weak form of Pleijel's theorem states that the number of eigenvalues $\lambda_j(\Omega)$, for which there exists an associated eigenfunction with precisely $j$ nodal domains (Courant-sharp eigenvalues), is finite. The purpose of this note is to determine an upper bound for Courant-sharp eigenvalues, expressed in terms of simple geometric invariants of $\Omega$. We will see that this is connected with one of the favorite problems considered by Y.~Safarov.
\end{abstract}\medskip

\noindent Keywords: {Dirichlet Laplacian, Nodal domains, Courant theorem,\\ Pleijel theorem.}\\[5pt]

%______________________________________________________

\section{Introduction and main result}\label{S-intro}

We consider the Dirichlet Laplacian $H(\Omega)$ in a bounded open set  $\Omega$ in $\mathbb R^d$. We denote by $\lambda_j(\Omega)$ ($j\in \mathbb N^*$) the sequence of eigenvalues arranged in nondecreasing order, with multiplicities. The ground state energy $\lambda_1(\Omega)$ is simply denoted by $\lambda(\Omega)$. We denote by $N(\phi_j) = \phi_j^{-1}(0)$ the \emph{nodal set} of an eigenfunction $\phi_j$ associated with $\lambda_j(\Omega)$, and by $\mu(\phi_j)$ the number of connected components of $\Omega \setminus N(\phi_j)$ (\emph{nodal domains)}.\\[5pt]
Courant's nodal domain theorem \cite[1923]{Cou} says that for any $j\ge 1$, the number  $\mu(\phi_j)$ is not greater than $j$. \medskip

 An eigenvalue $\lambda_j(\Omega)$ is called \emph{Courant-sharp} if there exists an associated eigenfunction $\phi_j$ with $\mu(\phi_j) = j$. In contrast with Sturm's theorem in dimension~$1$, the weak form of Pleijel's theorem \cite[1956]{Pl} says:

\begin{theorem}\label{Pleijelweakform}
In dimension $2$, the number of Courant-sharp eigenvalues of $H(\Omega)$ is finite.
\end{theorem}

This theorem is the consequence of a more precise theorem (strong Pleijel's theorem):

\begin{theorem}\label{PropPleijel1}
In dimension $2$, for any sequence of spectral pairs $(\phi_n,\lambda_n)$  of $H(\Omega)$,
\begin{equation}\label{compa}
\limsup_{n\rightarrow  +\infty} \frac{\mu(\phi_n)}{n}\leq \frac{4 \pi }{ \lambda (\mathbb D_1)} =\left(\frac{2}{{\bf j}} \right)^2 < 1\,,
\end{equation}
where $\mathbb D_1$ is the disk of unit area, and ${\bf j}$ the least positive zero of the Bessel function $J_0$.
\end{theorem}

\noindent \textbf{Remark}. Pleijel's theorem extends to bounded domains in $\mathbb{R}^d$, and more generally to compact $d$-dimensional manifolds with boundary, see Peetre \cite{Pe}, B\'{e}rard and Meyer \cite{BeMe}. More precisely for $d\ge 2$, there exists a constant $\gamma(d) < 1$ such that
\begin{equation}\label{compad}
\limsup_{n\rightarrow  +\infty} \frac{\mu(\phi_n)}{n}\leq \gamma(d)\,.
\end{equation}
It is interesting to note that the constant $\gamma(d)$ only depends on the dimension and is otherwise independent of the geometry.

In view of Pleijel's theorem, it is a natural question to look for geometric upper bounds for Courant-sharp eigenvalues. The purpose of this note is to give a \emph{geometrically controlled} version of Theorem~\ref{Pleijelweakform}. In dimension $2$, we prove the following result.

\begin{theorem}\label{thbound}
Let $\Omega \subset \mathbb{R}^2$ be a bounded $C^2$ domain. Then, there exists a positive constant $\beta(\Omega)$ depending only on the geometry of $\Omega$, such that any Courant-sharp eigenvalue $\lambda_k(\Omega)$ of $H(\Omega)$ satisfies
$$
k \frac{\lambda(\mathbb{D}_1)}{|\Omega|} \le \lambda_k(\Omega)  \leq  \beta(\Omega).
$$
More precisely, the constant $\beta(\Omega)$ can be computed in terms of the area $|\Omega|$, the peri\-me\-ter $\ell(\partial \Omega)$ of $\Omega$,  as well as bounds on the curvature of $\partial \Omega$ and on  the cut-distance\footnote{The cut-distance is defined in Section~\ref{S-sa}, Equation~\eqref{cut}.} $\varepsilon_0(\Omega)$   to $\partial \Omega$.
\end{theorem}%
The result also holds with weaker regularity assumptions. For example, inspection of the proof which uses the results of van den Berg and  Lianantonakis \cite{BL} gives the following non optimal but more explicit corollary.
\begin{corollary}\label{C-eb}
Let $\Omega \subset \mathbb{R}^2$ be a bounded domain. Define the geometric quantity $D(\Omega)$ by
$$
D(\Omega) = \sup_{\varepsilon > 0}\frac{\left| \left\lbrace  x \in \Omega ~:~ d(x) < \varepsilon \right\rbrace \right|}{\varepsilon} \,,
$$
where $d(x)$ is the distance from $x$ to the boundary of $\Omega$. If $D(\Omega)$ is finite, then any Courant-sharp eigenvalue $\lambda_k(\Omega)$ satisfies,
$$
\lambda(\D_1) \, k \le \,  |\Omega| \, \lambda_k(\Omega) \le 2 \, \left( \frac{24\,\pi\, \lambda(\D_1)}{\lambda(\D_1)-4\pi}\right)^4 \, \frac{\left(D(\Omega) \right)^4}{|\Omega|^2} \,.
$$
\end{corollary}
Observe that the lower and upper bounds are dilation invariant. When $\Omega$ is regular, $D(\Omega)$ can be bounded from above by
$$
D(\Omega) \le \max \left\lbrace \frac{|\Omega|}{\varepsilon_0(\Omega)}, 2 \ell(\partial \Omega) \right\rbrace .
$$

\noindent \textbf{Remarks}. (i) Corollary~\ref{C-eb} holds as soon as the boundary of $\Omega$ has Minkowski dimension $1$, see Section~\ref{S-BL}. (ii) The constant $D(\Omega)$ is bigger than the upper Minkowski content. We cannot substitute $D(\Omega)$ with the upper Minkowski content because we need  upper bounds on the quantities involved, not only an asymptotic behaviour.\medskip

In all the paper, we only consider the Dirichlet problem.  It would also be interesting to analyze the Neumann problem in the same spirit. Looking at the proof of Polterovich in \cite{Pol}, the main point would be to  obtain a geometric estimate of the number of nodal domains touching the boundary.\\

\noindent {\bf Organization of the paper.}\\
The paper is organized as follows. In Section~\ref{S-ppt}, we sketch the proofs of Pleijel's theorem, and we explain  the  idea of how to obtain geometric upper bounds for Courant-sharp eigenvalues. In Section~\ref{S-lbcf} we describe lower bounds on the counting function, using \cite{Saf} or \cite{BL}, and we derive upper bounds for the Courant-sharp eigenvalues. In Section~\ref{S-ex}, we compare the bounds obtained in Section~\ref{S-lbcf} for three very simple examples (the disk, the annulus and the square), and the bounds one can derive for other explicit examples (rectangles, equilateral triangles, etc.).\medskip

\noindent{\bf Added in proof.} We point out the following recent paper: M. van den Berg, K. Gittins, On the number of Courant sharp Dirichlet eigenvalues, arXiv 1602.08376. \medskip

\section{Proofs of Pleijel's theorem}\label{S-ppt}

In this section, we sketch the proof of Theorem~\ref{PropPleijel1} for a domain $\Omega$ in $\R^d$. We first introduce some notation.\\
Let $N_{\Omega}(\lambda)$ denote the counting function for $H(\Omega)$,
\begin{equation}\label{E-cf}
N_{\Omega}(\lambda) = \# \left\lbrace j ~|~ \lambda_j(\Omega) < \lambda \right\rbrace \,.
\end{equation}
The counting function can be written as
\begin{equation}\label{E-cfr}
N_{\Omega}(\lambda) = C_d \,  |\Omega|\,  \lambda^{\frac d2} - R(\lambda)\,,
\end{equation}
where $C_d$ is the Weyl constant, $|\Omega|$ denotes the $d$-dimensional volume of $\Omega$, and the remainder term $R(\lambda)$ satisfies $R(\lambda) = o(\lambda^{\frac d2})$ according to Weyl's theorem. The Weyl constant is given by
\begin{equation}\label{wconst}
C_{d}:= (2\pi)^{-d}\omega_d\,,
\end{equation}
where $\omega_d$ is the volume of the unit ball in $\R^d$,
\begin{equation}\label{eq:vol}
\omega_d = \pi^\frac d2 / \Gamma(\frac d2 +1)\,.
\end{equation}
We also denote by $\mathbb{B}_1^{d}$ the ball of volume $1$ in $\R^d$.\medskip

To prove Theorem~\ref{PropPleijel1}, we  start with the identity
\begin{equation}\label{idpl}
\frac{\mu(\phi_n)}{n}\, \frac{n}{\lambda_n(\Omega)^{\frac d2}} \, \frac{\lambda_n(\Omega)^{\frac d2}}{\mu(\phi_n)} =1\,.
\end{equation}
Applying the Faber-Krahn inequality to each nodal domain of $\phi_n$ and summing up, we have
\begin{equation} \label{FK}
\frac{\lambda_n(\Omega)^{\frac d2}}{\mu(\phi_n)}\geq \frac{ \lambda(\mathbb B_1^d)^{\frac d2}}{|\Omega|}\,.
\end{equation}
Note for  later reference that
\begin{equation} \label{FKa}
\text{if~} \mu(\phi_n)=n\,,\, \text{then~} \frac{\lambda_n(\Omega)^{\frac d2}}{n}\geq \frac{ \lambda(\mathbb B_1^d)^{\frac d2}}{|\Omega|}\,.
\end{equation}
This  gives a necessary condition for $\lambda_n(\Omega)$ to be Courant-sharp, which  is (up to the renormalization by the volume) independent of the geometry of $\Omega$.\\[3pt]
Taking a subsequence $\phi_{n_i}$ such that $$\lim_{i\rightarrow +\infty}\frac{\mu(\phi_{n_i})}{n_i} = \limsup_{n\rightarrow  +\infty} \frac{\mu(\phi_n)}{n}\,,$$
 and implementing in \eqref{idpl}, we deduce:
\begin{equation}\label{ineqpl1}
 \frac{ \lambda(\mathbb B_1^d)^{\frac d2}}{|\Omega|}\, \limsup_{n\rightarrow +\infty} \frac{\mu(\phi_n)}{n}\, \lim_{n\rightarrow + \infty} \frac{n}{N_{\Omega}(\lambda_n)}\, \lim_{\lambda \rightarrow +\infty} \frac{N_{\Omega}(\lambda)}{\lambda^{\frac d2}} \leq 1\,.
\end{equation}
Having in mind Weyl's formula, we obtain
\begin{equation}\label{E-gamma}
\limsup_{n\rightarrow +\infty} \frac{\mu(\phi_n)}{n} \le
\gamma(d) := \frac{1}{C_d \, \lambda(\mathbb B_1^d)^{\frac d2}}\,.
\end{equation}
When $d=2$, one has $C_2 = \frac{1}{4\pi}$, $\lambda(\mathbb B_1^2)=\pi {\bf j}^2$, so that $\gamma(2) = \frac{4}{{\bf j}^2} < 1$ since ${\bf j} \approx 2.40\,$. More generally, for $d\ge 2\,$, one has
$$
\gamma (d):=\frac{2^{d-2} d^2 \Gamma (d/2)^2}{(j_{\frac{d-2}{2},1})^d}\,,
$$
where $j_{\nu,1}$ denotes the first positive zero of the Bessel function $J_\nu$ (in particular $j_{0,1}={\bf j}$), and it can be shown, see \cite{BeMe}, that
\begin{equation}\label{E-fk}
\gamma(d) < 1\,.
\end{equation}
This proves Theorem~\ref{PropPleijel1}, and Theorem~\ref{Pleijelweakform} follows as well.\hfill $\square$ \medskip

\noindent \textbf{Remark}. In the case of general Riemannian manifolds, one needs to use an adapted isoperimetric inequality which is valid for domains with small enough  volume, see \cite{Pe, BeMe}.\medskip

We now give an alternative proof of Theorem~\ref{Pleijelweakform} which provides a hint on how to bound the Courant-sharp eigenvalues from above.\\[5pt]
If $\lambda_n$ is Courant-sharp, then $\lambda_{n-1} < \lambda_n$ and hence, $n = N_{\Omega}(\lambda_n) + 1$. Using \eqref{FKa}, we obtain
\begin{equation}\label{CF-4}
\lambda_n \text{~Courant-sharp~} \Rightarrow N_{\Omega}(\lambda_n) + 1 = n \le |\Omega| \, \left( \frac{\lambda_n(\Omega)}{\lambda(\mathbb B_1^d)}\right)^{\frac d2}\,.
\end{equation}
Writing the counting function as,
\begin{equation}\label{CF-6}
N_{\Omega}(\lambda) = C_d\, |\Omega|\, \lambda^{\frac d2} - R(\lambda)\,.
\end{equation}
and plugging this relation into \eqref{CF-4}, we obtain that
\begin{equation}\label{CF-8}
\lambda_n(\Omega) \text{~Courant-sharp~} \Rightarrow
F_{\Omega}\left( \lambda_n(\Omega)\right)\le 0\,,
\end{equation}
where the function $F_{\Omega}$ is defined for $\lambda > 0$ by
\begin{equation}\label{E-nc2}
F_{\Omega}(\lambda) = C_d \left( 1 - \gamma(d) \right)\, |\Omega|\, \lambda^{\frac d2} -R(\lambda) +1\,.
\end{equation}
By Weyl's theorem, the remainder term satisfies $ R(\lambda) = o(\lambda^{\frac{d}{2}})$.  Since $1 - \gamma(d) > 0$, see \eqref{E-fk}, the function $F_{\Omega}$ tends to infinity when $\lambda$ tends to infinity and hence the number of Courant-sharp eigenvalues must be finite.\hfill $\square$\medskip

As a matter of fact, the preceding proof tells us that Courant-sharp eigenvalues must be less than or equal to
\begin{equation}\label{E-nc4}
\inf \{\mu > 0 ~|~ F_{\Omega}(\lambda) > 0 \text{~for~} \lambda \ge \mu\} \,.
\end{equation}
Although this quantity is a geometric invariant associated with $\Omega$, it is not clear how to estimate it in terms of simple geometric invariants, even if we used Ivrii's sharp estimate  $R(\lambda) = \mathcal{O}(\lambda^{\frac{d-1}{2}})$, \cite{Iv, Se}. In order to proceed, it is sufficient to have an \emph{explicit geometric} upper bound  $\overline{R}(\lambda)$ of $R(\lambda)$. Indeed, define the function
\begin{equation}\label{E-nc6}
\overline{F}_{\Omega}(\lambda) = C_d \left( 1 - \gamma(d) \right) \, |\Omega| \, \lambda^{\frac d2} - \overline{R}(\lambda) +1\,.
\end{equation}
Then, any Courant-sharp eigenvalue $\lambda_k(\Omega)$ must satisfy $\overline{F}_{\Omega}\left( \lambda_k(\Omega) \right) \le 0$, and hence the inequality
\begin{equation}\label{E-nc8}
\lambda_k(\Omega) \le \inf \{\mu > 0 ~|~ \overline{F}_{\Omega}(\lambda) > 0 \text{~for~} \lambda \ge \mu\} \,.
\end{equation}

In the next section, we use the explicit upper bounds $\overline{R}(\lambda)$ provided by the papers of Safarov \cite{Saf} and van den Berg and Lianantonakis \cite{BL} to obtain upper bounds on the Courant-sharp eigenvalues in terms of simple geometric invariants.

\section{Lower bounds on the counting function and applications to Courant-sharp eigenvalues}\label{S-lbcf}

In this section, we describe lower bounds on the counting functions derived from \cite{Saf} or \cite{BL}, and apply them to bounding the Courant-sharp eigenvalues.

\subsection{The approach via Y. Safarov}\label{S-sa}

Here, we implement a result by Y. Safarov  \cite[2001]{Saf} which provides a lower bound for the spectral function on the diagonal, with an explicit control on the remainder term. This estimate is obtained by making use of finite propagation speed for the wave equation, and precise Tauberian theorems.\\
If $\Omega \subset \mathbb R^d$ is an open set, then the spectral function of the Dirichlet Laplacian
$$
e(x,x,\lambda):= \frac{1}{2} \left( \sum_{\lambda_j < \lambda} \phi_j(x)^2  +  \sum_{\lambda_j \leq \lambda} \phi_j(x)^2\right) \,,
$$
satisfies \cite[Cor. 3.1]{Saf}
\begin{equation}\label{saflb}
e(x,x,\lambda) \geq C_d \lambda^\frac d 2 - \frac{2 d\, C_d \pi ^{-1}\,\nu_{m_d}^2}{d(x)} \left(\lambda^\frac 12 + \frac{\nu_{m_d}}{d(x)}\right)^{d-1}\,,
\end{equation}
for all $x \in \Omega$ and $\lambda > 0\,$.\\
Here  $d(x)$ is the Euclidean distance to $\partial \Omega$,
and  $\nu_{m_d}$ is a universal constant depending only on the dimension.\\
More precisely,  let
$$
m_d = \left\{\begin{array}{ll}\frac{d+1}{2}\,,& \mbox{ if } d \mbox{ is odd},\\[5pt]
\frac{(d+2)}{2} \,,&  \mbox{ if } d \mbox{ is even}.
\end{array}
\right.
$$
Then,
$$
\nu_m = (\widetilde \nu_m)^{\frac {1}{2m}}\,,
$$
where $\widetilde \nu_m$ is the ground state energy of the Dirichlet realization of $(-1)^m \frac{d^{2m}}{dt^{2m}}$ on $]-\frac 12,\frac 12[$.\medskip

Define
$$
\widetilde N_{\Omega}(\lambda):= \int_\Omega e(x,x,\lambda) \,dx \,,
$$
and let $ \epsilon_0(\Omega)$   be the largest  number $\varepsilon$ with the property that $$\Omega_{\epsilon}^{b}:= \{x\in \Omega\,,\, d(x) < \epsilon\}$$ is diffeomorphic to $ \partial \Omega \times ]0,\epsilon [$.\\
 Then, for $0 < \epsilon <\epsilon_0(\Omega)$,
$$
\begin{array}{ll}
\widetilde N_{\Omega}(\lambda) \geq & C_d \,  |\Omega|\,\lambda^\frac d 2 - C_d \, |\Omega_\epsilon^b|   \lambda^\frac d 2 \\
& - 2 d\, C_d \,  \pi ^{-1} \nu_{m_d}^2 \,  \left(\lambda^\frac 12 + \frac{\nu_{m_d}}{\epsilon}\right)^{d-1}\,\left(\int _{d(x) >\epsilon} \frac 1 {d(x)} dx  \right)\,.
\end{array}
$$

\noindent This inequality is also true by semi-continuity for $N_{\Omega}(\lambda)$.\medskip

Writing $N_{\Omega}(\lambda) = C_d \,  |\Omega|\,\lambda^\frac d 2 - R(\lambda)$ as in \eqref{CF-6}, we have
\begin{equation}\label{RS}
R(\lambda) \le C_d \, |\Omega_\epsilon^b|   \lambda^\frac d 2 + 2 d\, C_d \,  \pi ^{-1} \nu_{m_d}^2 \,  \left(\lambda^\frac 12 + \frac{\nu_{m_d}}{\epsilon}\right)^{d-1}\,\left(\int _{d(x) >\epsilon} \frac 1 {d(x)} dx  \right)\,.
\end{equation}

We now use our freedom for choosing $\epsilon$. A  convenient choice in order to get the right power of $\lambda$ is to take
\begin{equation}
\epsilon := \alpha(\Omega)   \lambda^{-\frac 12}\,.
\end{equation}
Because we need this estimate for any $\lambda$ in the spectrum of the Laplacian,  and actually for $\lambda > \lambda_2(\Omega)$ (because the Courant-sharp property is already established for the two first eigenvalues), we choose
\begin{equation}\label{opti}
 \alpha(\Omega) = \epsilon_0(\Omega) \lambda_2(\Omega)^\frac 12\,.
\end{equation}
To have more explicit bounds, we could also choose
$$
\alpha (\Omega)=\epsilon_0(\Omega) \underline{\lambda_2}(\Omega)^\frac 12 \,
$$
where $\underline{\lambda_2}(\Omega)$ is a geometric lower bound of
$\lambda_2 (\Omega)$  (using Faber-Krahn inequality or a consequence of Li-Yau inequality, see below \eqref{eq:FK} and  \eqref{eq:LY}).\\[5pt]
For regular domains, the right-hand side of \eqref{RS} can be estimated in terms of the geometry of $\Omega$. \medskip

\noindent \textbf{For the sake of simplicity, we give the details in the case $d=2$}.\medskip

\noindent In dimension $2$, the above lower bound for $\widetilde{N}_{\Omega}(\lambda)$ (and $N_{\Omega}(\lambda)$) reads
\begin{equation}\label{D2-N0}
\begin{array}{ll}
N_{\Omega}(\lambda) & \geq C_2 |\Omega|   \lambda  - C_2  \left| \Omega_{\alpha \lambda^{-\frac 12}}^{b}\right| \lambda\\
 &   - 4 C_2 \pi^{-1} \nu_{2}^2 \,  \left( 1 + \frac{\nu_{2}}{\alpha}\right) \,  \lambda^\frac 12\, \left(\int_{\left\{ d(x) >\alpha \lambda^{-\frac 12}\right\} } \, \frac 1 {d(x)} dx \right)\,.
\end{array}
\end{equation}
When $d=2$, we have $m_2= 2$, and we can verify (using the quasimode $(\frac 14 -x^2)^2$ of \cite{Saf})  that
$$\widetilde \nu_2 \leq  7\times 8 \times 9 \leq 2^9\,,$$
which implies  the rough estimate
$$ \nu_2 \leq 4 \cdot \, 2^{\frac 14}\leq 5 \,.$$

We now assume  that $\partial \Omega$ is a smooth submanifold, so that $\partial \Omega$ is the union of $p$ smooth simple closed curves. We write the proof in the case $p=1$, the general case is similar. Let $c : [0,L] \to \mathbb{R}^2$ be a parametrization of $\partial \Omega$ by arc-length, with $L := \ell(\partial \Omega)$, the length of the boundary. The associated Frenet frame is $\{\tau(s), \nu(s)\}$. We can assume that the orientation is chosen such that $\nu(s)$ points towards the interior of $\Omega$. The curvature $\kappa(s)$ of the curve is given by the equation $\dot{\tau}(s) = \kappa(s) \nu(s)$. Let $\kappa_{-}(\Omega)$ denote the infimum of $\kappa$ over $[0,L]$. Define the map
$$\left\{
\begin{array}{l}
F: [0,L] \times ]-\infty, \infty[ \to \R^2 \,,\\
F(s,t) = c(s) + t \nu(s)\,.
\end{array}
\right.
$$

We have
$$\partial_s F(s,t) \wedge \partial_t F(s,t) = \left( 1 - t \kappa(s)\right) \tau(s) \wedge \nu(s)\,.
$$

\noindent The map $F$ is a local diffeomorphism for $|t| < t_+$  with
$$ t_+:= \left( \sup_{[0,L]}|\kappa(s)| \right)^{-1}\,.$$
The injectivity of $F$ is determined by the infimum $ \underline{\delta}_{+}$ of the \emph{cut-distance} $ \delta_{+}(s)$ to the submanifold $\partial \Omega$, where
\begin{equation}\label{cut}
\delta_+(s) := \sup\{t > 0 ~:~ t = \mathrm{dist}(F(s,t),\partial \Omega)\}.
\end{equation}
 In this case, we have
\begin{equation}
 \epsilon_0(\Omega) = \inf\{t_{+},\underline{\delta}_{+}\}
 \end{equation}
 so that $F$ is a diffeomorphism from $[0,L]\times ]0,\epsilon_0(\Omega)[$ onto its image (i.e. so that $F$ is both a local diffeomorphism and injective). For $\epsilon < \epsilon_0(\Omega)$, we have
$$
|\Omega_{\epsilon}^{b}| = \int_{0}^{L} \int_{0}^{\epsilon} \left( 1 - t \kappa(s) \right) \, ds \,dt\,.
$$
It follows that

\begin{equation}\label{D2-N2}
\left\{
\begin{array}{l}
  C_2 \left| \Omega_{\alpha \lambda^{\frac 12}}^{b} \right| \, \lambda \le \beta_1(\Omega) \lambda^{\frac 12}\,, \text{where}\\[5pt]
\beta_1(\Omega) = \frac{1}{4\pi} \left( 1+\epsilon_0(\Omega) |\kappa_{-}(\Omega)|\right) \epsilon_0(\Omega) \lambda_2(\Omega)^{\frac 12} \ell (\partial \Omega)\,.
\end{array}
\right.
\end{equation}
The third term in the right-hand side of \eqref{D2-N0} can be written as

\begin{equation}\label{D2-N3}
\left\{
\begin{array}{l}
\beta_2(\Omega) \lambda^{\frac 12} \int_{\left\{ d(x) >\alpha \lambda^{-\frac 12}\right\} } \, \frac 1 {d(x)} dx\,, \text{where}\\[7pt]
\beta_2(\Omega) := \pi^{-2} \nu_2^2 \left( 1+ \nu_2 \epsilon_0(\Omega)^{-1}\lambda_2(\Omega)^{- \frac 12}\right)\,.
\end{array}
\right.
\end{equation}
Write
$$
\begin{array}{ll}
\int_{\left\{ d(x) >\alpha \lambda^{-\frac 12}\right\} } \, \frac 1 {d(x)} dx = &\int_{\left\{ d(x) > \epsilon_0(\Omega)\right\} } \, \frac 1 {d(x)} dx \\&+ \int_{0}^{L} \int_{\alpha \lambda^{-\frac 12}}^{\epsilon_0(\Omega)} \frac{\left( 1-t \kappa(s)\right)}{t}ds\,dt
\end{array}
$$
We can estimate the second integral in the right-hand side as we did  above. The first integral can be estimated from above by $|\Omega|/\varepsilon_0(\Omega)$. It follows that there exist positive constants $\beta_3(\Omega)$
and $\beta_4(\Omega)$ such that, for all $\lambda > \lambda_2(\Omega)$,
$N_{\Omega} (\lambda) = \frac{|\Omega|}{4\pi}\lambda - R(\lambda)$, with

\begin{equation}\label{D2-N6}
R(\lambda) \le \overline{R}(\lambda) = \beta_3(\Omega) \lambda^\frac 12  \ln \left( \frac{\lambda}{\lambda_2(\Omega)}\right) + \beta_4(\Omega) \lambda^{\frac 12}\,.
\end{equation}
Note that the constants only depend on the geometry of the domain $\Omega$.  More precisely, the constants can be computed in terms of $|\Omega|, \ell(\partial \Omega)$, $\kappa_{-}(\Omega), \epsilon_0(\Omega)$, and $\lambda_2(\Omega)$.\medskip

\noindent \textbf{Remarks}.  (i) The preceding proof shows that one can alternatively estimate the constants in terms of $|\Omega|, \ell(\partial \Omega)$, $\epsilon_0(\Omega)$, $\lambda_2(\Omega)$, and the number of holes of the domain (through the integral $\int_{\partial \Omega} \kappa$).\\(ii) In higher  dimensions, one can state a similar result in which the curvature $\kappa$ of the curve is replaced by the mean curvature $h$ of the hypersurface $\partial \Omega$. For this purpose, one uses the Heintze-Karcher comparison theorem \cite{HK78}.\medskip

\noindent Applying \eqref{E-nc6} and \eqref{E-nc8}, we obtain that any Courant-sharp eigenvalue $\lambda_k(\Omega)$ satisfies
\begin{equation}\label{ineq}
f_{\Omega}(\lambda_k) \le 0\,,
\end{equation}
where the function $f_{\Omega}$ is defined for $\mu > \lambda_2(\Omega)$ by
\begin{equation}
f_{\Omega}(\mu) := \frac{\lambda(\mathbb D_1)- 4\pi}{4\pi \lambda(\mathbb D_1)} \,|\Omega| \mu - \beta_3(\Omega) \mu ^\frac 12 \ln\left(\frac{\mu}{\lambda_2}\right)  - \beta_4(\Omega) \mu^{\frac 12}+1 \,.
\end{equation}

Since $\lambda(\D_1) > 4\pi$, see  \eqref{E-fk},  the coefficient of the term $\mu$ in the expression of $f_{\Omega}$ is positive, so that the function tends to infinity when $\mu$ tends to infinity. Hence $I_\Omega:= f_\Omega^{-1} (]-\infty,0])$ is either empty or bounded from above.

Define
$$
\beta_S(\Omega)=\max\{\lambda_2(\Omega),\beta_0(\Omega)\}\,,
$$
where $\beta_0(\Omega)$ is  the supremum of  $I_\Omega $ if $I_\Omega$ is non empty and $0$ otherwise. From Equation \eqref{ineq} we conclude that
\begin{equation}\label{b-Sa}
\lambda_k(\Omega) \text{~Courant-sharp~} \Rightarrow
\lambda_k(\Omega) \le \beta_S(\Omega)\,.
\end{equation}
We have proved Theorem \ref{thbound}. \hfill $\square$\medskip

Starting from the inequality $f_{\Omega}(\lambda_k) \le 0$ in the above proof, we conclude that any Courant-sharp eigenvalue $\lambda_k$ satisfies
$$
A_2\, |\Omega| \lambda_k^{\frac 12} \le \beta_3(\Omega) \ln \frac{\lambda_k}{\lambda_2} + \beta_4(\Omega)\,,
$$
where $A_2 = \frac{1}{4\pi} - \frac{1}{\lambda(\D_1)}$. Using the inequality $\ln \frac{\mu}{\lambda_2} \le 4\,\left(\frac{\mu}{\lambda_2}\right)^{\frac 14}$ which holds for any $\mu \ge \lambda_2$, we obtain the following more explicit bound.

\begin{corollary}\label{C-Sa}
In dimension $2$, any Courant-sharp eigenvalue $\lambda_k(\Omega)$ of $H(\Omega)$ satisfies
\begin{equation}\label{sa-exp}
\lambda_k(\Omega) \le \max \left\lbrace \lambda_2(\Omega), \left( \frac{16\pi \lambda(\D_1)}{\lambda(\D_1) - 4\pi} \right)^4 \, \frac{\left( \beta_3(\Omega)+\beta_4(\Omega)\right)^4}{|\Omega|^4\lambda_2(\Omega)}\right\rbrace\,.
\end{equation}
\end{corollary}%

\noindent \textbf{Remarks}. (i) For the unit disk, the bound \eqref{sa-exp} is sharper than Corollary~\ref{C-eb}, see Section~\ref{S-ex}. (ii) P\'olya's conjecture for Dirichlet eigenvalues (see \cite{Poly}) does not go in the right direction. Indeed lower bounds on the Dirichlet eigenvalues correspond to upper bounds on $N(\lambda)$. This would be good for Neumann eigenvalues, but in this case, there are other problems, see \cite{Pl} and more recently \cite{Pol}.\medskip
%
%In the next Section, we obtain an estimate similar to \eqref{b-Sa} by %using a lower bound of the counting function due to van den %Berg--Lianantonakis, \cite{BL}.

\subsection{Approach via van den Berg--Lianantonakis}\label{S-BL}

Prior to Y. Safarov, van den Berg and Lianantonakis have given lower bounds for the counting function $N_{\Omega}(\lambda)$ depending on the
Minkowski dimension of $\partial \Omega$. When this dimension is $(d-1)$, they prove \cite[Theorem 2.1]{BL} that if
\begin{equation}\label{MvdB1}
\lambda \geq 4 |\Omega|^{-\frac 2 d} \,,
\end{equation}
then
\begin{equation}\label{MvdB2}
N (\lambda) \geq C_d |\Omega| \, \lambda^\frac d2
- 3 D(\Omega) \lambda^{(d-1)/2} \log \left((2|\Omega|)^{\frac 2 d}
\, \lambda \right)\,,
\end{equation}
where the geometric constant $D(\Omega)$ is defined by
\begin{equation}\label{defD}
D(\Omega):= \sup_\epsilon \frac{ |\Omega_\epsilon^b|}{\epsilon}\,.
\end{equation}

To apply \eqref{MvdB2} to Pleijel's theorem, one needs to compare condition \eqref{MvdB1} with the condition $\lambda > \lambda_2(\Omega)$.
One can for example observe that the  Faber-Krahn  inequality applied  to the second eigenvalue gives (see \cite{AsBe} or \eqref{FKa} for $d=2$)
\begin{equation}\label{eq:FK}
\lambda_2(\Omega) \geq  (2 \omega_d)^{\frac 2 d}   |\Omega|^{-\frac 2 d} j_{\frac d2 -1,1}^2\,.
\end{equation}
For $d=2$, since $2 \pi j_{0,1}^2 >4$, the condition $\lambda >\lambda_2(\Omega)$ implies \eqref{MvdB1}. For $d\geq 2$,  we can use the following lower bound for $\lambda_2(\Omega)$ which is a consequence of Li-Yau inequality (see   \cite[Formula (11.5)]{AsBe}),
\begin{equation}\label{eq:LY}
\lambda_2(\Omega) > \frac{d}{d+2} \frac{4 \pi^2 2^{\frac 2 d}}{ (\omega_d
|\Omega|)^{\frac 2 d}}\,.
\end{equation}

Hence it is enough to verify that:
$$
\frac{d}{d+2} \frac{4 \pi^2 2^{\frac 2 d}}{ (\omega_d)^{\frac 2 d}} \geq 4\,,
$$
which is easy to establish. Indeed, using \eqref{eq:vol}, we obtain
$$
\frac{d}{d+2}\, \pi \, 2^{\frac 2 d}\, \Gamma(\frac d2 +1) ^{\frac 2 d} \geq 1\,,
$$
which follows from the inequality $\frac{d}{d+2}\, \pi \geq 1$ for $d\geq 1$.\\

\noindent \textbf{Assuming $d=2$ for the sake of simplicity}, and using \eqref{CF-8} together with \eqref{MvdB2}, we obtain that any Courant-sharp eigenvalue $\lambda_k(\Omega)$, with $\lambda_k > \lambda_2$, satisfies $g_{\Omega}(\lambda_k) \le 0$, where $g_{\Omega}$ is defined by
$$
g_{\Omega}(\mu) = \left( \frac{1}{4\pi} - \frac{1}{\lambda(\D_1)}\right) |\Omega|\,  \lambda - 3 D(\Omega)\, \mu^{\frac 12}\ln(2|\Omega|\mu)+1\,,
$$
for $\mu \ge \lambda_2(\Omega)$. Define $\beta_1(\Omega)$ to be $0$ if $g_{\Omega}(\mu) \ge 0$, and $\sup\{\mu > \lambda_2 ~:~ g_{\Omega}(\mu) \le 0\}$ otherwise, and define
$$
\beta_B(\Omega) := \max \{\lambda_2(\Omega),\beta_1(\Omega)\}\,.
$$
Then,
\begin{equation}\label{b-BL}
\lambda_k(\Omega) \text{~Courant-sharp~} \Rightarrow \lambda_k(\Omega) \le \beta_B(\Omega)\,.
\end{equation}
This proves Theorem~\ref{thbound} using the lower bound for the counting function provided by \cite{BL}.\hfill $\square$\medskip

From the inequality $g_{\Omega}(\lambda_k) \le 0$ in the preceding proof, we have that any Courant-sharp eigenvalue $\lambda_k(\Omega)$ satisfies the inequality
$$
\left(\frac{1}{4 \pi}  - \frac{1}{\lambda(\D_1)} \right) |\Omega|\, \lambda_k^\frac 12 - 3D(\Omega) \ln (2 |\Omega|\lambda_k) \leq 0\,,
$$
for $\lambda_k \ge \lambda_2$. Using the inequality
$$
\ln \mu \le 2 \, \mu^{\frac 14} \text{~~for~~} \mu \ge 16\,,
$$
and the fact that $2|\Omega|\lambda_k > 16$ (Faber-Krahn), we obtain  the more explicit bound given in Corollary~\ref{C-eb}.

%\begin{corollary}\label{C-expl}
%In dimension $2$, any Courant-sharp eigenvalue $\lambda_k(\Omega)$
%of $H(\Omega)$ satisfies the dilation-invariant inequality
%%
%\begin{equation}\label{E-sb}
%|\Omega|\lambda_k(\Omega) \le 2\, \left( \frac{24\pi %\lambda(\D_1)}{\lambda(\D_1)-4\pi}\right)^4 \, \frac{\left( D(\Omega) %\right)^4}{|\Omega|^2} \,.
%\end{equation}
%
%\end{corollary}%
%%
%To conclude, it remains to analyze $D(\Omega)$.

\noindent As kindly communicated by M.~van den Berg, in dimension $2$, when $\Omega$ is sufficiently regular, the  geometric invariant $D(\Omega)$ can be bounded from above by
$$
D(\Omega) \leq \max \left( \frac{|\Omega|}{\epsilon_0(\Omega)}\,,\, \ell (\partial \Omega) + \pi \, \varepsilon_0(\Omega) \, h(\Omega)  \right),
$$
where $h(\Omega)$ is the number of holes of $\Omega$, or by
$$
D(\Omega) \leq \max \left( \frac{|\Omega|}{\epsilon_0(\Omega)}\,,\, 2 \ell (\partial \Omega)  \right).
$$

%----------------------

\section{Examples and particular cases}\label{S-ex}

\subsection{Examples}

In some $2$-dimensional cases, it is possible to compute the upper bounds for Courant-sharp eigenvalues arising from the preceding sections explicitly. Consider the following domains,
$$
\begin{array}{l}
\Omega_1 = B(0,1), \text{~~the unit disc in~} \R^2\,,\\
\Omega_2 = B(0,1) \setminus B(0,a), 0<a<1, \text{~~the annulus~} A(0,a,1) \subset \R^2\,,\\
\Omega_3 = ]0,\pi[\times ]0,\pi[\, , \text{~~the square in~} \R^2 \text{~with side~} \pi\,.
\end{array}%
$$
For the unit disc, one finds that $\beta_S(\Omega_1) \approx 7.1\cdot 10^6$ and $\beta_B(\Omega_1) \approx 2.1\cdot 10^7$.\\%[5pt]
For the annulus, one finds that $\beta_B(\Omega_2) \approx 4.2\cdot 10^8$ when $a = 0.75\,$, \break and $\beta_B(\Omega_2) \approx 4\cdot 10^7$ when $a = 0.25\,$.  This indicates that the cut-distance to the boundary does matter in the upper bound on Courant-sharp eigenvalues.\\%[5pt]
For the square with side $\pi$, one finds that $\beta_B(\Omega_3) \approx 5.9\cdot  10^6$. It turns out that this bound is much bigger than the bound which is deduced in the next sub-section, namely $51$.\\
This is not surprising. The general lower bounds for the counting functions used in the preceding sections, Equations~\eqref{D2-N6} and \eqref{MvdB2}, are worse than the sharp $2$-dimensional estimate $R(\lambda) = \mathcal{O}(\lambda^{\frac 12})$, see \cite{Iv}, by a $\ln(\lambda)$ factor. On the other-hand, the estimate \eqref{Nlambdarec} has the right powers, and almost the right second constant.\\
Generally speaking one should therefore expect that the bounds $\beta_S(\Omega)$ and \break $\beta_B(\Omega)$ are not sharp.\medskip

\subsection{Particular cases}

As already mentioned,  improved Weyl's formulas with control of the remainder which are only asymptotic are not sufficient for an explicit version of Pleijel's theorem.
  We nevertheless mention for comparison  a formula due to V. Ivrii in 1980 (cf \cite[Chapter XXIX, Theorem 29.3.3 and Corollary 29.3.4]{Hor4}) which reads:
\begin{equation}\label{eq.Weyl2terms}
N(\lambda) = \frac{\omega_d}{(2\pi)^d} |\Omega|\, \lambda^\frac d 2 - \frac 14 \frac{\omega_{d-1}}{(2\pi)^{d-1}} \, |\partial \Omega|\,  \lambda^{\frac{d-1}{2}} + r(\lambda),
\end{equation}
where $r(\lambda)= \mathcal O (\lambda^{\frac{d-1}{2}} )$ in general, but can also be shown to be $o (\lambda^{\frac{d-1}{2}} )$ if the boundary is $C^\infty$, and under some generic conditions on the geodesic billiards (the measure of periodic trajectories should be zero). For piecewise smooth boundaries, see \cite{Va}. The second term is meaningful in this case only. \\

Formula \eqref{eq.Weyl2terms} is also established for irrational rectangles as a very special case in \cite{Iv}, but more explicitly in \cite{Ku} without any assumption of irrationality.  See also  \cite{Be83} for some $2$-dimensional domains with negative curvature. We do not discuss here the case of ``rough" boundaries which was in particular analyzed  by Netrusov et Safarov in \cite{NeSa} (and references therein).\\

Note that when $d=2$, the second term in \eqref{eq.Weyl2terms}  is
\begin{equation}\label{w2}
W_2(\lambda):=-\frac{1}{4\pi}|\partial \Omega|\, \lambda^{\frac{1}{2}} \,.
\end{equation}

 The Dirichlet (and Neumann) eigenvalues are explicitly given for few domains. In dimension $2$ these domains include the rectangles, the right-angled isosceles triangle, the equilateral triangle and the hemiequilateral triangle. In these cases, estimating the counting function amounts to estimating the number of points with integer coordinates inside some ellipse (these domains are obtained as quotient of a torus). The estimates which are obtained in this manner are compatible with Weyl's two terms asymptotic formula \eqref{eq.Weyl2terms}), involving the area of the domain and the length of it boundary. Similarly, in higher dimensions, one can explicitly describe the Dirichlet (and Neumann) eigenvalues of the fundamental domains of crystallographic affine Weyl groups, \cite{Be80}.  As far as the asymptotic estimate is concerned, this is possible because the remainder term in Weyl's estimate has order $ \lambda^{\frac{d-2}{2}+ \frac{1}{d+1}}$ for a $d$-dimensional torus. \medskip

\paragraph{\bf Rectangle.}~\\
Following (and improving) a remark in a course of R. Laugesen \cite{BeHe1}, one has a lower bound of $N(\lambda)$ in the case of the rectangle $\mathcal R=\mathcal R(a,b):=(0,a\pi) \times (0,b\pi)$, which can be expressed in terms of area and perimeter. One can indeed observe that the area of the intersection of the ellipse $\{ \frac{(x+1)^2}{a^2} + \frac{(y+1)^2}{b^2} <\lambda\}$ with $\mathbb R^+\times \mathbb R^+$ is a lower bound for $N(\lambda)$.\\
 The formula reads:
\begin{equation}\label{Nlambdarec}
N_{\mathcal R}(\lambda) > \frac{1}{4\pi} |\mathcal R|\, \lambda - \frac{1}{2\pi}|\partial\mathcal R |\,  \sqrt{\lambda} +1\,, \text{~~for~}\lambda \geq \frac{1}{a^2}+\frac{1}{b^2}\,.
\end{equation}
Here we can observe that the second term is $2\, W_2(\lambda)$ (see \eqref{w2}).

\paragraph{\bf Equilateral triangle, see \cite{BeHe4}.}~\\
We consider the equilateral triangle with side $1$.
\begin{equation}\label{ETR-12}
N_{\mathcal T}(\lambda) \ge \frac{\sqrt{3}}{4} \frac{\lambda}{4\pi} - \frac{3}{2\pi} \sqrt{\lambda} + 1\,.
\end{equation}
Again we  observe that the second term is $2\, W_2(\lambda)$ (see \eqref{w2}).

\paragraph{\bf Right-angled isosceles triangle, see \cite{BeHe4}.}~\\
Let $\mathcal B_{\pi}$ denote the right-angled isosceles triangle,
\begin{equation}\label{OTB-4}
\mathcal B_{\pi} = \left\lbrace (x,y) \in ]0,\pi[^2 ~|~ y < x \right\rbrace\,.
\end{equation}
\begin{equation}\label{OTB-10}
N_{\mathcal B}(\lambda) \ge \frac{\pi \lambda}{8} - \frac{(4+\sqrt{2})\sqrt{\lambda}}{4} - \frac{1}{2}\,.
\end{equation}

\paragraph{\bf The cube, see \cite{HK}.}~\\
For the cube $]0,\pi[^3$, we have,  for $\lambda \ge 3$\,:
\begin{equation}\label{lauzz}
N(\lambda) >   \frac \pi 6 \lambda^\frac 32 - \frac{3\pi}{4} \lambda + 3  \sqrt{\lambda -2} -1\,.
\end{equation}

\noindent{\bf Acknowledgements.\\}
The authors would like to thank M. van den Berg for useful discussions and remarks on a preliminary version of this paper, and the referee for his comments.\\[5pt]
%----------------------

\bibliographystyle{plain}

\end{document}